\documentclass[11pt,bezier]{article}
\usepackage{amsmath,amssymb,amsfonts, euscript}
\newpage
\newcommand{\be}{\begin{equation}}
\newcommand{\ee}{\end{equation}}
\newtheorem{theorem}{Theorem}[section]
\newtheorem{lemma}[theorem]{Lemma}
\newtheorem{corollary}{Corollary}[section]
\newtheorem{definition}{Definition}[section]
\newtheorem{remark}{Remark}[section]
\newtheorem{example}{Example}[section]
\newtheorem{proposition}[theorem]{Proposition}

\renewcommand{\theequation}{\arabic{section}.\arabic{subsection}.\arabic{equation}}
\title{\bf\Large Spectral Analysis of Multi-dimensional Self-similar Markov Processes }
\vspace{-1cm}
\author
{ N. Modarresi\,\,\,\, and \,\,\,  S. Rezakhah\thanks{ Faculty of
Mathematics and Computer Science, Amirkabir University of
Technology, 424 Hafez Avenue, Tehran 15914, Iran. E-mail: $\;$
namomath@aut.ac.ir (N. Modarresi),$\;\;$ rezakhah@aut.ac.ir (S.
Rezakhah).} }
\date{}
\begin{document}
\maketitle

\begin{abstract}
In this paper we consider a discrete scale invariant (DSI) process $\{X(t), t\in {\bf R^+}\}$ with  scale $l>1$. We consider to have some fix number of observations in every scale, say $T$, and to get our samples at discrete points $\alpha^k,\; k\in {\bf W}$ where $\alpha$ is obtained by the equality $l=\alpha^T$ and ${\bf W}=\{0, 1, \ldots\}$. So we provide a discrete time scale invariant (DT-SI) process $X(\cdot)$ with parameter space $\{\alpha^k, k\in {\bf W}\}$.
We find the spectral representation of the covariance function of such DT-SI process. By providing harmonic like representation of multi-dimensional self-similar processes, spectral density function of them are presented.
We assume that the process $\{X(t), t\in {\bf R^+}\}$ is also Markov in the wide sense and provide a discrete time scale invariant Markov (DT-SIM) process with the above scheme of sampling. We present an example of DT-SIM process, simple Brownian motion, by the above sampling scheme and verify our results.
Finally we find the spectral density matrix of such DT-SIM process and show that its associated $T$-dimensional self-similar Markov process is fully specified by $\{R_{j}^H(1),R_{j}^H(0),j=0, 1, \ldots, T-1\}$ where $R_j^H(\tau)$ is the covariance function of $j$th and $(j+\tau)$th observations of the process.\\


\end{abstract}

\section{Introduction}
The concept of stationarity and self-similarity are used as a fundamental property to handle many natural phenomena. Lamperti transformation defines a one to one correspondence between stationary and self-similar processes. A function is scale invariant if it is identical to any of its rescaled version, up to some suitable renormalization in amplitude.
Discrete scale invariance (DSI) process can be defined as the Lamperti transform of periodically correlated (PC) process \cite {b2}. Many critical systems, like statistical physics, textures in geophysics, network traffic and image processing can be interpreted by these processes \cite {b1}. Flandrin et. al. introduced a multiplicative spectral representation of DSI processes based on the Mellin transform and presented preliminary remarks about estimation issues \cite{b2}, \cite{f1}.

As the Fourier transform is known as a suited representation for stationarity, but not for self-similarity. A harmonic like representation of self-similar process is introduced by using Mellin transform \cite {f1}.

Covariance function and spectral density of discrete time periodically correlated Markov process have been
studied and characterized in \cite {n1}. Markov processes have been the center of extensive research activities and wide sense Markov processes are studied before the general theory. In some texts, these processes are defined in the case of transition probabilities of a Markov process. Various classes of wide sense Markov processes are like jump processes, diffusion processes and processes with a discrete interference of chance \cite {g0}.
A process which is Markov and self-similar, is called self-similar Markov process. These processes are involved in various parts of probability theory, such as branching processes and fragmentation theory \cite {c1}.

In this paper, we consider a DSI process with some scale $l>1$, and we get our samples at points $\alpha^k$, where $k\in {\bf W}$, $l=\alpha^T$, ${\bf W}=\{0, 1, \ldots\}$ and $T$ is the number of samples in each scale. By such sampling we provide a discrete time scale invariant process in the wide sense and find the spectral representation of the covariance function of such process.

This paper is organized as follows.
In section 2, we present stationary and self-similar processes by shift and renormalized dilation operators. Then we provide a suitable platform for our study of discrete time self-similar (DT-SS) and discrete time scale invariant (DT-SI) processes by introducing quasi Lamperti transformation. Harmonizable representation of these processes are expressed in this section too.
Also by using the spectral density matrix of PC processes, the spectral representation of the covariance function of DT-SI processes are given.
In section 3, a harmonic like representation of multi-dimensional self-similar processes and spectral density function of them are obtained. As an example we introduce a process called simple Brownian motion which is DSI and Markov too.
Finally a discrete time scale invariant Markov (DT-SIM) process with the above scheme of sampling is considered in section 3 and the spectral density matrix of such process and its associated $T$-dimensional self-similar Markov process are characterized.

\vspace{3mm}
\renewcommand{\theequation}{\arabic{section}.\arabic{equation}}
\section{Theoretical framework}
\setcounter{equation}{0}
In this section, by using renormalized dilation operator, we define discrete time self-similar and discrete time scale invariant processes. The quasi Lamperti transformation and its properties are introduced. We also present the harmonizable representation of stationary and harmonic like representation of self-similar processes. The spectral density of PC processes and the spectral representation of the covariance function of DT-SI processes are given.

\subsection{Stationary and self-similar processes}

\begin{definition}
Given $\tau\in {\bf R}$, the shift operator ${\cal S}_{\tau}$ operates on process $\{Y(t),t\in {\bf R}\}$ according to \be({\cal S}_{\tau}Y)(t):=Y(t+\tau).\ee
A process $\{Y(t),t\in {\bf R}\}$ is said to be stationary, if for any $t,\tau\in {\bf R}$
\be \{({\cal S}_{\tau}Y)(t)\}\stackrel{d}{=}\{Y(t)\}\ee
where $\stackrel{d}{=}$ is the equality of all finite-dimensional distributions.\\ \\
If $(2.2)$ holds for some $\tau \in {\bf R}$, the process is said to be periodically correlated. The smallest of such $\tau$ is called period of the process.
\end{definition}

\begin{definition}
Given some numbers $H>0$ and $\lambda>0$, the renormalized dilation operator ${\cal D}_{H,\lambda}$ operates on process $\{X(t),t\in {\bf R^+}\}$ according to
\be({\cal D}_{H,\lambda}X)(t):=\lambda^{-H}X(\lambda t).\ee
A process $\{X(t),t\in {\bf R^+}\}$ is said to be self-similar of index $H$, if for any ${\lambda}>0$
\be \{({\cal D}_{H,\lambda}X)(t)\}\stackrel{d}{=}\{X(t)\}.\ee
The process is said to be {\em DSI} of index $H$ and scaling factor ${\lambda}_0>0$ or {\em (H,${\lambda}_0$)-DSI}, if $(2.4)$ holds for $\lambda=\lambda_0$.
\end{definition}

As an intuition, self-similarity refers to an invariance with respect to any dilation factor. However, this may be a too strong requirement for capturing in situations that scaling properties are only observed for some preferred dilation factors.\\

\begin{definition}
A process $\{X(k),k\in {\bf \check{T}}\}$ is called discrete time self-similar {\em (DT-SS)} process with parameter space $\check{T}$,
where $\check{T}$ is any subset of distinct points of positive real numbers, if for any $k_1, k_2 \in \check{T}$
\be \{X(k_2)\}\stackrel{d}{=}(\frac{k_2}{k_1})^H\{X(k_1)\}.\ee
The process $X(\cdot)$ is called discrete time scale invariance {\em (DT-SI)} with scale $l>0$ and parameter space $\check{T}$, if for any $k_1, k_2=lk_1 \in \check{T}$, $(2.5)$ holds.
\end{definition}

\begin{remark}
If the process $\{X(t),t\in {\bf R^+}\}$ is {\em DSI} with scale $l=\alpha^T$ for fixed $T\in {\bf N}$ and $\alpha>1$, then by sampling of the process at points $\alpha^{k}, k\in {\bf W}$ where ${\bf W}=\{0, 1, \ldots\}$, we have $X(\cdot)$ as a {\em DT-SI} process with parameter space $\check{T}=\{\alpha^{k}; k\in {\bf W}\}$ and scale $l=\alpha^T$.
If we consider sampling of $X(\cdot)$ at points $\alpha^{nT+k}, n\in {\bf W}, \mbox{for fixed}\ k=0, 1,\ldots, T-1$, then $X(\cdot)$ is a {\em DT-SS} process with parameter space $\check{T}=\{\alpha^{nT+k}; n\in {\bf W}\}$.
\end{remark}

Yazici et.al. \cite {y1} and \cite {w1} introduced wide sense self-similar processes as the following definition, which can be obtained by applying the Lamperti transformation ${\cal L}_{H}$ to the class of wide-sense stationary processes. This class  encompasses all strictly self-similar processes with finite variance, including Gaussian processes such as fractional Brownian motion but no other alpha-stable processes.

\begin{definition}
A random process $\{X(t),t\in {\bf R^+}\}$ is said to be wide sense self-similar with index H, for some $H>0$ if the following properties are satisfied for each $a>0$\\

$(i)\,\,\ E[X^2(t)]<\infty$,

$(ii)\,\,E[X(at)]=a^HE[X(t)]$,

$(iii)\,\, E[X(at_1)X(at_2)]=a^{2H}E[X(t_1)X(t_2)]$.\\ \\
This process is called wide sense {\em DSI} of index $H$ and scaling factor $a_0>0$, if the above conditions hold for some $a=a_0$.\\
\end{definition}

\begin{definition}
A random process $\{X(k),k\in \check{T}\}$ is called {\em DT-SS} in the wide sense with index $H>0$ and with parameter space $\check{T}$, where $\check{T}$ is any subset of distinct points of positive real numbers, if for all $k, k_1\in \check{T}$ and all $a>0$, where $ak, ak_1\in \check{T}:$\\

$(i)\,\,\ E[X^2(k)]<\infty$,

$(ii)\,\,E[X(ak)]=a^HE[X(k)]$,

$(iii)\,\, E[X(ak)X(ak_1)]=a^{2H}E[X(k)X(k_1)]$.\\ \\
If the above conditions hold for some fixed $a=a_0$, then the process is called {\em DT-SI} in the wide sense with scale $a_0$.
\end{definition}

\begin{remark}
Let $\{X(t),t\in {\bf R^+}\}$ in Remark $2.1$ be {\em DSI} in the wide sense. Then $X(\cdot)$ with parameter space $\check{T}=\{\alpha^{k}; k\in {\bf W}\}$ for $\alpha>1$ is {\em DT-SI} in the wide sense, where ${\bf W}=\{0, 1, \ldots\}$ and $X(\cdot)$ with parameter space $\check{T}=\{\alpha^{nT+k}; n\in {\bf W}\}$ for fixed $T\in {\bf N}$ and $\alpha>1$ is {\em DT-SS} in the wide sense.
\end{remark}

Through this paper we are dealt with wide sense self-similar and wide sense scale invariant process, and for simplicity we omit the term "in the wide sense" hereafter.

\subsection{Quasi Lamperti transformation}
We introduce the quasi Lamperti transformation and its properties by followings.
\begin{definition}
The quasi Lamperti transform with positive index $H$ and $\alpha>1$, denoted by ${\cal L}_{H,\alpha}$ operates on a random process
$\{Y(t),t\in {\bf R}\}$ as
\be {\cal L}_{H,\alpha}Y(t)=t^HY(\log_\alpha t)\ee
and the corresponding inverse quasi Lamperti transform ${\cal L}^{-1}_{H,\alpha}$ on process $\{X(t), t\in {\bf R^+}\}$
acts as
\be{\cal L}^{-1}_{H,\alpha}X(t)={\alpha}^{-tH}X(\alpha^t).\ee
\end{definition}
\vspace{5mm}
One can easily verify that ${\cal L}_{H,\alpha}{\cal L}^{-1}_{H,\alpha}X(t)=X(t)$ and
${\cal L}^{-1}_{H,\alpha}{\cal L}_{H,\alpha}Y(t)=Y(t).$ Note that in the above definition, if $\alpha=e$, we have the usual Lamperti transformation ${\cal L}_{H}$.

\begin{theorem}
The quasi Lamperti transform guarantees an equivalence between the shift operator ${\cal S}_{\log_{\alpha}k}$ and the renormalized dilation operator ${\cal D}_{H,k}$ in the sense that, for any $k>0$
\be {\cal L}^{-1}_{H,\alpha}{\cal D}_{H,k}{\cal L}_{H,\alpha}={\cal S}_{\log_{\alpha}k}.\ee
\end{theorem}
{\bf Proof:}
$${\cal L}^{-1}_{H,\alpha}{\cal D}_{H,k}{\cal L}_{H,\alpha}Y(t)={\cal L}^{-1}_{H,\alpha}{\cal D}_{H,k}(t^HY(\log_{\alpha}t))={\cal L}^{-1}_{H,\alpha}(k^{-H}(kt)^HY(\log_{\alpha}{kt}))$$
$$={\cal L}^{-1}_{H,\alpha}(t^HY(\log_{\alpha}{kt}))={\alpha}^{-tH}({\alpha}^{t})^HY(\log_{\alpha}{k{\alpha}^{t}})=
Y(\log_{\alpha}k+t)={\cal S}_{\log_{\alpha}k}Y(t).\square$$

\vspace{3mm}
\begin{corollary}
If $\{Y(t),t\in {\bf R}\}$ is stationary process, its quasi Lamperti transform  $\{{\cal L}_{H,\alpha}Y(t),\\t\in {\bf R^+}\}$ is self-similar. Conversely if $\{X(t),t\in {\bf R^+}\}$ is self-similar process, its inverse quasi Lamperti transform $\{{\cal L}^{-1}_{H,\alpha}X(t),t\in {\bf R}\}$ is stationary.\\
\end{corollary}

\begin{corollary}
If $\{X(t), t\in {\bf R^+}\}$ is {\em $(H,{\alpha}^T$)-DSI} then ${\cal L}^{-1}_{H,\alpha}X(t)=Y(t)$ is {\em PC} with period $T>0$. Conversely if $\{Y(t), t\in {\bf R}\}$ is {\em PC} with period $T$ then ${\cal L}_{H,\alpha}Y(t)=X(t)$ is {\em $(H,{\alpha}^T$)-DSI}.\\
\end{corollary}

\begin{remark}
If $X(\cdot)$ is a {\em DT-SS} process with parameter space $\check{T}=\{l^n, n\in {\bf W}\}$, then its stationary counterpart $Y(\cdot)$ has parameter space $\check{T}=\{nT, n\in {\bf N}\}$:
$$X(l^n)={\cal L}_{H,\alpha}Y(l^n)=l^{nH}Y(\log_{\alpha}{\alpha^{nT}})=\alpha^{nTH}Y(nT).$$\\
Also it is clear by the following relation that if $X(\cdot)$ is a {\em DT-SI} process with scale $l=\alpha^T$, $T\in {\bf N}$ and parameter space $\check{T}=\{\alpha^k, k\in {\bf W}\}$, then $Y(\cdot)$ is a discrete time periodically correlated {\em (DT-PC)} process with period $T$ and parameter space $\check{T}=\{n, n\in {\bf N}\}$:
$$Y(n)={\cal L}^{-1}_{H,\alpha}X(n)=\alpha^{-nH}X(\alpha^n).$$
\end{remark}

\subsection{Harmonizable representation}
A stationary process $Y(t)$, $EY(t)=0$, can be represented as
\be Y(t)=\int_{-\infty}^{\infty}e^{i\omega t}d\varphi(\omega)\ee
which is called harmonizable representation of the process, and $\varphi(\omega)$ is a right continuous orthogonal increment process, see \cite {l1}. Also the covariance function can be represented as
\be R_Y(t,s)=\int_{-\infty}^{\infty}\int_{-\infty}^{\infty}e^{it\omega-is\omega'}d\Phi(\omega,\omega')\ee
where the spectral measure satisfies
\be d\Phi(\omega,\omega')=E[d\varphi(\omega)\overline{d\varphi(\omega')}]=\left\{\begin {array}{cc}
\hspace{3mm}0\hspace{12mm}\omega\neq\omega'\\
d\Psi(\omega)\hspace{7mm}\omega=\omega'\\
\end {array}\right.\ee
and $d\Psi(\omega)=E[|d\varphi(\omega)|^2]$. All the spectral mass is located on the diagonal line $\omega=\omega'$. When $\Phi(\omega,\omega')$ is absolutely continuous, we have spectral density $\phi(\omega,\omega')$ such that $d\Phi(\omega,\omega')=\phi(\omega,\omega')d\omega d\omega'.$
A necessary and sufficient condition for this equality to hold, as Loeve's condition for harmonizability, is that $\Phi(\omega,\omega')$ must satisfy
$\int\int |d\Phi(\omega,\omega')|<\infty.$
The corresponding notion for processes after a Lamperti transformation introduces a new representation for a class of
processes deviating from self-similarity, which is called multiplicative harmonizability. A self-similar process $X(t)$ has harmonic like representation as an inverse Mellin transform, namely an integral of uncorrelated spectral increments $d\varphi(\omega)$ on the Mellin basis \cite {b1}.
\be X(t)=\int t^{H+i\omega}d\varphi(\omega)\ee
and the process has this property if it verifies as
\be R_X(t,s)=\int\int t^{H+i\omega}s^{H-i\omega'}d\Phi(\omega,\omega').\ee
The inverse Mellin transformation gives  the expression of the spectral function if the correlation is known as \cite {b2}
\be \phi(\omega,\omega')=\int\int t^{-H-i\omega}s^{-H+i\omega'}R_X(t,s)\frac{dt}{t}\frac{ds}{s}.\ee

\subsection{Spectral density function}
The spectral density of a PC process is introduced by Gladyshev in \cite {g1}. If $Y(n)$ is a DT-PC process, the spectral density matrix is Hermitian nonnegative definite $T\times T$ matrix of functions $f(\omega)=[f_{jk}(\omega)]_{j,k=0,1,\ldots,T-1}$, and
the covariance function has the representation
\be R_n(\tau):=\mathrm{Cov}\big(Y(n),Y(n+\tau)\big)=\sum_{k=0}^{T-1}B_k(\tau)e^{2k\pi in/T}\ee
\vspace{-1mm}
where
$$B_k(\tau)=\int_{0}^{2\pi}e^{i\tau\omega}f_k(\omega)d\omega.$$
Also $f_k(\omega)$ and $f_{jk}(\omega), j,k=0,1,\ldots,T-1$ are related through
$$f_{jk}(\omega)=\frac{1}{T}f_{k-j}\big((\omega-2\pi j)/T\big),\hspace{5mm}0\leqslant\omega<2\pi .$$
For $k<0$, $\omega<0$ or $\omega>2\pi$, the functions $f_k(\omega)$ are defined by the equality $f_{k+T}(\omega)=f_k(\omega)$ and $f_k(\omega)=f_k(\omega+2\pi)$.\\

Let $\{X(t), t\in {\bf R^+}\}$ be a zero mean DSI process with scale $l$. If $l<1$, we reduce the time scale, so that $l$ in the new time scale is greater than $1$. Our sampling scheme is to get samples at points $\alpha^k$, $k\in {\bf W}$, where by choosing the number of samples in each scale, say $T\in {\bf N}$, we find $\alpha$ by $l=\alpha^T$. Therefore the process under study $\{X(\alpha^n), n\in{\bf W}\}$ is DT-SI with scale $l=\alpha^T$.

\begin{proposition}
If $X(\alpha^n)$ is {\em DT-SI} with scale $l=\alpha^T$, $T\in {\bf N}$, then we have the spectral representation of the covariance function of the process as
\be R_n^H(\tau):=\mathrm{Cov}\big(X(\alpha^n),X(\alpha^{n+\tau})\big)=\alpha^{(2n+\tau)H}\sum_{k=0}^{T-1}B_k(\tau)e^{2k\pi in/T}\ee
where
\vspace{-5mm}
\be B_k(\tau)=\int_{0}^{2\pi}e^{i\tau\omega}f_k(\omega)d\omega\ee
and
\be f_{jk}(\omega)=\frac{1}{T}f_{k-j}\big((\omega-2\pi j)/T\big)\ee
for $j,k=0,1,\ldots,T-1$ and $0\leqslant\omega<2\pi$.\\
\end{proposition}
{\bf Proof:} According to (2.6) and Corollary 2.1, for any $n,\tau\in {\bf W}$
$$R_n^H(\tau)=E[X(\alpha^n)X(\alpha^{n+\tau})]=E[{\cal L}_{H,\alpha}Y(\alpha^n){\cal L}_{H,\alpha}
Y(\alpha^{n+\tau})]$$
$$=\alpha^{(2n+\tau)H}
E[Y(n)Y(n+\tau)]$$
where $Y(n)$ is DT-PC process with period $T=\log_{\alpha}l$. Thus by (3.1)
$$R_n^H(\tau)=\alpha^{(2n+\tau)H}R_n(\tau)=\alpha^{(2n+\tau)H}\sum_{k=0}^{T-1}B_k(\tau)e^{2k\pi in/T}.\square$$

\renewcommand{\theequation}{\arabic{section}.\arabic{equation}}
\section{Characterization of the spectrum}
\setcounter{equation}{0}
In this section we provide spectral density matrix of multi-dimensional self-similar process $W(n)$. By using harmonic like representation of a self-similar process, we characterize the spectral density matrix of DT-SI process in subsection 3.1.
A discrete time scale invariant Markov (DT-SIM) process with a new scheme of sampling is considered and the properties of an introduced example is verfied too. The spectral density matrix of such process and its associated $T$-dimensional self-similar
Markov process are characterized in subsection 3.2.

\subsection{Spectral representation of multi-dimensional self-similar process}
By Rozanov \cite {r1}, if $\xi(t)=\{\xi^k(t)\}_{k=1,\ldots,n}$ be an $n$-dimensional stationary process, then
\be \xi(t)=\int e^{i\lambda t}\phi(d\lambda)\ee
is its spectral representation, where $\phi=\{\varphi_k\}_{k=1,\ldots,n}$ and $\varphi_k$ is the random spectral measure associated with the $k$th component $\xi^k$ of the $n$-dimensional process $\xi$. Let
$$B_{kr}(\tau)=E[\xi^k(\tau+t)\overline{\xi^r(t)}], \hspace{7mm}k,r=1,\ldots,n$$
and $B(\tau)=[B_{kr}(\tau)]_{k,r=1,\ldots,n}$ be the correlation matrix of $\xi$.
The components of the correlation matrix of the process $\xi$ can be represented as
\be B_{kr}(\tau)=\int e^{i\lambda \tau}F_{kr}(d\lambda), \hspace{7mm}k,r=1,\ldots,n\ee
where for any Borel set $\Delta$, $F_{kr}(\Delta)=E[\varphi_k(\Delta)\overline{\varphi_r(\Delta)}]$  are the complex valued set functions which are $\sigma$-additive and have bounded variation. For any $k,r=1,\ldots,n$, if the sets $\Delta$ and $\Delta'$ do not intersect, $E[\varphi_k(\Delta)\overline{\varphi_r(\Delta')}]=0$.
For any interval $\Delta=(\lambda_1,\lambda_2)$ when $F_{kr}(\{\lambda_1\})=F_{kr}(\{\lambda_2\})=0$ the following relation holds
\be F_{kr}(\Delta)=\frac{1}{2\pi}\int_{\Delta}\sum_{\tau=-\infty}^{\infty}B_{kr}(\tau)e^{-i\lambda \tau}d\lambda\ee
$$=\frac{1}{2\pi}B_{kr}(0)[\lambda_2-\lambda_1]+\lim_{T\rightarrow\infty}\frac{1}{2\pi}\sum_{0<|\tau|\leqslant T}
B_{kr}(\tau)\frac{e^{-i\lambda_2 \tau}-e^{-i\lambda_1 \tau}}{-i\tau}$$
in the discrete parameter case, and
$$F_{kr}(\Delta)=\lim_{a\rightarrow\infty}\frac{1}{2\pi}\int_{-a}^{a}\frac{e^{-i\lambda_2 \tau}-e^{-i\lambda_1 \tau}}{-i\tau}B_{kr}(\tau)d\tau$$
in the continuous parameter case.\\

Using the above results of Rozanov for multi-dimensional stationary processes and using the Lamperti transformation,
we present the definition of multi-dimensional self-similar process and obtain the properties of the corresponding multi-dimensional self-similar process by the following theorem.

\begin{definition}
The process $U(t)=\big(U^0(t), U^1(t),\ldots, U^{q-1}(t)\big)$ is a q-dimensional discrete time self similar process in the wide sense with parameter space $\check{T}$, which consists of finite or countable many points of ${\bf R}^+$, if the followings are satisfied\\

$\bf (a)$\hspace{3mm} $\{U^j(\cdot)\}$ for every $j=0, 1, \cdots, q-1$ is {\em DT-SS} process with parameter space $\check{T}$.\\

$\bf (b)$\hspace{3mm} $U^i(\cdot)$ and $U^j(\cdot)$ for $i, j=0, 1, \cdots, q-1$ have self-similar correlation, that is\\
$$\mathrm{Cov}\big(U^i(ts), U^j(tr)\big)=t^{2H}\mathrm{Cov}\big(U^i(s),U^j(r)\big).$$
\hspace{1cm}where $s, r, ts, tr$ are in $\check{T}$.\\
\end{definition}

\begin{theorem}
Let $W(\alpha^k)=\big(W^0(\alpha^k), W^1(\alpha^k),\ldots, W^{q-1}(\alpha^k)\big)$, $k\in {\bf Z},\;\alpha >1\,$ be a discrete time  $q$-dimensional self-similar process. Then\\

{\em (i)} The harmonic like representation of $W^j(\alpha^k)$ is \be
W^j(\alpha^k)= \alpha^{kH}\int_{0}^{2\pi}e^{i\omega
k}d\varphi_j\omega). \ee
where $\varphi_j(\omega)$ is the corresponding spectral measure, that  $E[d\varphi_j(\omega)\overline{d\varphi_r(\omega')}]=dD^H_{jr}(\omega)$, $j, r=0, 1,\ldots, q-1$ when $\omega=\omega'$ and is $0$ when $\omega\neq\omega'$. We call $D^H_{jr}(\omega)$ the spectral distribution function of the process.\\

{\em (ii)} The corresponding spectral density matrix of
$\{W(\alpha^k), k\in {\bf Z}\}$ is ${\bf d}^H(\omega)=[{\bf d}_{jr}^H(\omega)]_{j,r= 0, \ldots,q-1}$, where
\be {\bf d}^H_{jr}(\omega)=\frac{1}{2\pi}\sum_{n=-\infty}^{\infty}\alpha^{-nH}e^{-i\omega n}Q_{jr}^H(\alpha^n)\ee
$\alpha>1$ and $Q_{jr}^H(\alpha^n)$ is the covariance function of $W^j(\alpha^n)$ and $W^r(1)$.
\end{theorem}

Before proceeding to the proof of the theorem we remind that based on our sampling scheme at points $\alpha^k$,
$k\in {\bf Z}$, of continuous DSI process with scale $l=\alpha^T$, $\alpha \in {\bf R}$, $T\in N $. So 
we consider $W(\cdot)$ at points $l^n=\alpha^{nT}$ as the corresponding $T$-dimensional DT-SS process and apply this theorem in Lemma 3.4.\\ \\
{\bf Proof of (i):} $W^j(\alpha^k)$ for $j=0, 1,\ldots, q-1$ is
DT-SS and its stationary counterpart $\xi^j(k)$ has spectral
representation $\xi^j(k)=\int_{0}^{2\pi}e^{i\omega
k}d\varphi_j(\omega)$. Thus by (2.6)
$$W^j(\alpha^k)={\cal L}_{H,\alpha}\xi^j(\alpha^k)=\alpha^{kH}\xi^j(k)=\alpha^{kH}\int_{0}^{2\pi}e^{i\omega k}
d\varphi_j(\omega)
.$$\\
{\bf Proof of (ii):} The covariance matrix is denoted by
$Q^H(n,\tau)=[Q_{jr}^H(\alpha^n,\alpha^\tau)]_{j,r= 0, \ldots, q-1}$ where its
elements have the spectral representation
\be Q_{jr}^H(\alpha^m,\alpha^\tau)=E[W^j(\alpha^m\alpha^\tau)\overline{W^r(\alpha^m
)}]=\alpha^{2mH}E[W^j(\alpha^\tau)\overline{W^r(1)}]=\alpha^{2mH}Q_{jr}^H(\alpha^\tau).\ee
Also by (3.4)
\be Q_{jr}^H(\alpha^\tau)=\alpha^{\tau
H}E[\int_{0}^{2\pi}e^{i\omega\tau}d\varphi_j(\omega)\int_{0}^{2\pi}
\overline{d\varphi_r(\omega')}]=\alpha^{\tau
H}\int_{0}^{2\pi}e^{i\omega \tau}dD^H_{jr}(\omega)\ee
where $E[d\varphi_j(\omega)\overline{d\varphi_r(\omega')}]=dD^H_{jr}(\omega)$
when $\omega=\omega'$ and is $0$ when $\omega\neq\omega'$.\\

The spectral distribution function of the correlation matrix $Q^H(\alpha^k)=[Q_{jr}^H(\alpha^k)]_{j,r=
0, \ldots, q-1}$ is
$$D^H(d\omega)=[D_{jr}^H(d\omega)]_{j,r= 0, \ldots, q-1}.$$
(3.2)-(3.3), (3.7) and appropriate transformation we have
\be D_{jr}^H(A)=\frac{1}{2\pi}\int_{A}\sum_{n=-\infty}^{\infty}\alpha^{-nH}e^{-i\lambda n}Q_{jr}^H(\alpha^n)d\lambda.\ee
Let $A=(\omega,\omega+d\omega]$, then we have the spectral density matrix as
$d^H(\omega)=[d_{jr}^H(\omega)]_{j,r= 0, \ldots, q-1}$ where
$$d^H_{jr}(\omega):=\frac{D_{jr}^H(d\omega)}{d\omega}=\frac{1}{2\pi}\sum_{n=-\infty}^{\infty}\big(\frac{1}{d\omega}
\int_{\omega}^{\omega+d\omega}e^{-i\lambda
n}d\lambda\big)\alpha^{-Hn}Q_{jr}^H(\alpha^n)$$
$$=\frac{1}{2\pi}\sum_{n=-\infty}^{\infty}\big(-\frac{1}{in}\lim_{d\omega\rightarrow 0}
\frac{e^{-i{(\omega+d\omega)n}}-e^{-i\omega
n}}{d\omega}\big)\alpha^{-Hn}Q_{jr}^H(\alpha^n)$$
$$=\frac{1}{2\pi}\sum_{n=-\infty}^{\infty}\big((-\frac{1}{in})(-in)
e^{-i\omega
n}\big)\alpha^{-Hn}Q_{jr}^H(\alpha^n)=\frac{1}{2\pi}\sum_{n=-\infty}^{\infty}\alpha^{-nH}e^{-i\omega
n}Q_{jr}^H(\alpha^n).$$
Existence of ${\bf d}^H_{jr}(\omega)$ follows from part (i) of the theorem as $W^k(\alpha^n)$ is the Lamperti counterpart of the stationary process $\xi^k(n)$, $k=0, \ldots, q-1$.$\square$

\subsection{Spectral density of DT-SIM process}
Let $\{X(t), t\in {\bf R}\}$ be a DSI process with scale $l$ and Markov in the wide sense.
Using our sampling scheme described in this section, we assume $l$ and $\alpha$ to be greater than one.
Thus $\{X(\alpha^n), n\in{\bf W}\}$ is a discrete time scale invariant Markov (DT-SIM) process with scale $l=\alpha^T$.

Let $R(t_1,t_2)$ be some function defined on $\mathcal{T}\times\mathcal{T}$ and suppose that $R(t_1,t_2)\neq 0$ everywhere on $\mathcal{T}\times\mathcal{T}$, where $\mathcal{T}$ is an interval. Borisov \cite {b3} showed that the necessary and sufficient condition for $R(t_1,t_2)$ to be the covariance function of a Gaussian Markov process with time space $\mathcal{T}$ is that
\be R(t_1,t_2)=G\big(\min(t_1,t_2)\big)K\big(\max(t_1,t_2)\big)\ee
where $G$ and $K$ are defined uniquely up to a constant multiple and the ratio $G/K$ is a positive nondecreasing function on $\mathcal{T}$.

It should be noted that the Borisov result on Gaussian Markov processes can be easily derived in the discrete case for second order Markov processes in the wide sense, by using Theorem 8.1 of Doob \cite {d1}.\\

Here we present a closed formula for the covariance function of the DT-SIM process and characterized the covariance matrix of corresponding $T$-dimensional self-similar Markov process by theorems 3.2 and 3.3 \cite{m1}.
\begin{theorem}
Let $\{X(\alpha^n),n\in {\bf Z}\}$ be a {\em DT-SIM} process with
scale $l=\alpha^T$, $\alpha>1$, $T\in {\bf N}$, then the covariance
function \be
R_n^H(\tau)=E[X(\alpha^{n+\tau})X(\alpha^n)]\ee\vspace{3mm} where
$\tau\in{\bf Z}$, $n=0, 1, \ldots, T-1$,
$R_{n+T}^H(\tau)=\alpha^{2TH}R_n^H(\tau)$ and $R_n^H(\tau)\neq0$ is
of the form \be
R_n^H(kT+v)=[\tilde{h}(\alpha^{T-1})]^k\tilde{h}(\alpha^{v+n-1})[\tilde{h}(\alpha^{n-1})]^{-1}R_n^H(0)\ee
$$R_n^H(-kT+v)=\alpha^{-2kTH}R_{n+v}^H((k-1)T+T-v)$$
where $k\in{\bf Z}$, $v=0,1,\ldots,T-1$,
\be\tilde{h}(\alpha^r)=\prod_{j=0}^{r}h(\alpha^j)=\prod_{j=0}^{r}R_j^H(1)/R_j^H(0),\hspace{7mm}r\in{\bf Z}\ee
\vspace{-2mm}
and $\tilde{h}(\alpha^{-1})=1$.\\
\end{theorem}
{\bf Proof:} Here we present the sketch of the proof. From the Markov property (3.9), for $\alpha>1$, we have that\\
\be R_n^H(\tau)=G(\alpha^n)K(\alpha^{n+\tau})\hspace{1cm}\tau\in {\bf W}.\ee
and
$$R_{n}^H(0)=G(\alpha^n)K(\alpha^n).$$
Thus
\be K(\alpha^{n+\tau})=\frac{R_{n}^H(\tau)}{R_{n}^H(0)}K(\alpha^n).\ee
For $\tau=1$, by a recursive substitution in (3.14) one can easily verify that
\be K(\alpha^n)=K(1)\prod_{j=0}^{n-1}h(\alpha^j)\ee
where $h(\alpha^j)=R_{j}^H(1)/R_{j}^H(0).$
Hence for $n=0,1,\ldots,T-1,\hspace{2mm}k\in{\bf W}$
$$K(\alpha^{kT+n})=K(1)\prod_{j=0}^{kT+n-1}h(\alpha^j).$$
As $X(\cdot)$ is DT-SI with scale $\alpha^T$ by (3.10)
$$h(\alpha^{T+i})=\frac{R_{T+i}^H(1)}{R_{T+i}^H(0)}=\frac{R_{i}^H(1)}{R_{i}^H(0)}=h(\alpha^i),\hspace{1cm}i\in {\bf W}.$$
Therefore using (3.12) we have
\be\prod_{j=0}^{kT+n-1}h(\alpha^j)=[\tilde{h}(\alpha^{T-1})]^k\tilde{h}(\alpha^{n-1}).\ee
Consequently for $n=0,1,\ldots,T-1$
\be K(\alpha^{kT+n})=K(1)\bigg[\tilde{h}(\alpha^{T-1})\bigg]^k\tilde{h}(\alpha^{n-1}).\ee\\
Let $\tau=kT+v$, then it follows from (3.14) and (3.17) that
$$R_n^H(kT+v)=\frac{K(\alpha^{n+kT+v})}{K(\alpha^n)}R_{n}^H(0)=\frac{K(1)[\tilde{h}(\alpha^{T-1})]^k
\tilde{h}(\alpha^{v+n-1})}{K(1)\tilde{h}(\alpha^{n-1})}R_{n}^H(0)$$
$$=[\tilde{h}(\alpha^{T-1})]^k\tilde{h}(\alpha^{v+n-1})[\tilde{h}(\alpha^{n-1})]^{-1}R_n^H(0)$$
for $k=0, 1,\ldots$, $\alpha>1$ and $n,v=0, 1,\ldots, T-1$. Similar to the above proof for $\tau=-kT+v$ we have
$$R_n^H(-kT+v)=[\tilde{h}(\alpha^{T-1})]^{-k}\tilde{h}(\alpha^{v+n-1})[\tilde{h}(\alpha^{n-1})]^{-1}R_n^H(0)$$
and also note that
$$R_n^H(-kT+v)=E[X(\alpha^{-kT+n+v})X(\alpha^n)]=\alpha^{-2kTH}E[X(\alpha^{n+v})X(\alpha^{kT+n})]$$
$$=\alpha^{-2kTH}R_{n+v}^H(kT+v)=\alpha^{-2kTH}R_{n+v}^H((k-1)T+T-v).\square$$

\begin{example}
We consider moving of a particle in different environment $A_1, A_2, \ldots$ based on Brownian motion with different rates. Specially we consider this movement by $X(t)$ with index $H>0$ and scale $\lambda>1$ as
$$X(t)=\sum_{n=1}^{\infty}\lambda^{n(H-\frac{1}{2})}I_{[\lambda^{n-1}, \lambda^n)}(t)B(t)$$
where $B(\cdot)$, $I(\cdot)$ are Brownian motion and indicator function respectively and we call this process, simple Brownian motion.  \end{example}

Let $A_1=[1, \lambda)$, $A_2=[\lambda, \lambda^2)$ and $A_n=[\lambda^{n-1}, \lambda^n)$ as disjoint sets.
The process $X(t)$ is DSI and Markov too. For checking these properties, first we find the covariance function of it.
The covariance function of the process for $t\in A_n$, $s\in A_m$ and $s\leqslant t$ is
\be\mathrm{Cov}\big(X(t), X(s)\big)=\lambda^{(n+m)(H-\frac{1}{2})}\mathrm{Cov}\big(B(t), B(s)\big)=\lambda^{(n+m)(H-\frac{1}{2})}s\ee
since as we know $\mathrm{Cov}\big(B(t), B(s)\big)=\min\{t,s\}$. Therefore by the condition $(3.9)$, the above covariance is the covariance function of a Markov process. Now we verify the DSI property. If $t$ is in $(\lambda^{n-1}, \lambda^n]$ then $\lambda t$ is in $(\lambda^{n}, \lambda^{n+1}]$. Thus for $t\in A_{n+1}$ and $s\in A_{m+1}$ we have
\vspace{3mm}
$$\mathrm{Cov}\big(X(\lambda t), X(\lambda s)\big)=\lambda^{(n+m+2)(H-\frac{1}{2})}\mathrm{Cov}\big(B(\lambda t), B(\lambda s)\big)=\lambda^{(n+m+2)(H-\frac{1}{2})}\lambda s$$
$$=\lambda^{2H}\lambda^{(n+m)(H-\frac{1}{2})}s=\lambda^{2H}\mathrm{Cov}\big(X(t), X(s)\big).$$
Then $X(t)$ is (H,${\lambda}$)-DSI.\\

By sampling of the process $X(\cdot)$ at points $\alpha^n$, $n\in {\bf W}$, where $\lambda=\alpha^T$, $T\in {\bf N}$ and $\lambda>1$, we provide a DT-SIM process and investigate the conditions of Theorem $3.2$.
For $j=kT+i$ where $i=0, 1, \ldots, T-2$ and $k=0, 1, \ldots$ by $(3.18)$ we have that
$$h(\alpha^j)=\frac{R_j^H(1)}{R_j^H(0)}=\frac{\mathrm{Cov}\big(X(\alpha^{j+1}), X(\alpha^j)\big)}{\mathrm{Cov}\big(X(\alpha^j), X(\alpha^j)\big)}=\frac{\alpha^{2(k+1)TH'+j}}{\alpha^{2(k+1)TH'+j}}=1$$
as $\alpha^j, \alpha^{j+1} \in A_{k+1}$ and  $H'=H-\frac{1}{2}$. Also for $j=kT+T-1$ we have that
$$h(\alpha^j)=\frac{R_j^H(1)}{R_j^H(0)}=\frac{\mathrm{Cov}\big(X(\alpha^{j+1}), X(\alpha^j)\big)}{\mathrm{Cov}\big(X(\alpha^j), X(\alpha^j)\big)}=\frac{\alpha^{(2k+3)TH'+j}}{\alpha^{(2k+2)TH'+j}}=\alpha^{TH'}$$
as $\alpha^j\in A_{k+1}$ and $\alpha^{j+1}\in A_{k+2}$. Thus for $j=kT+i$, $i=0, 1, \ldots, T-2$ and $k=0, 1, \ldots$
$$\tilde{h}(\alpha^{kT+i})=\prod_{r=0}^{kT+i}h(\alpha^r)=\prod_{r=0}^{k}\alpha^{TH'}=\alpha^{kTH'}$$
and for $j=kT+T-1$
$$\tilde{h}(\alpha^{kT+T-1})=\prod_{r=0}^{kT+T-1}h(\alpha^r)=\prod_{r=0}^{k+1}\alpha^{TH'}=\alpha^{(k+1)TH'}.$$
Finally as $\tilde{h}(\alpha^{T-1})=\alpha^{TH'}$,
$$\tilde{h}(\alpha^{v+n-1})=\left\{\begin {array}{cc}
1\hspace{2cm}v+n-1\leqslant T-2\\
\alpha^{TH'}\hspace{15mm}v+n-1\geqslant T-1\\
\end {array}\right.$$
and $\tilde{h}(\alpha^{n-1})=1$, $R_n^H(0)=E[X(\alpha^n)X(\alpha^n)]=\alpha^{2TH'+n}$, $n=0, 1, \ldots, T-1$. Thus\\
$$R_n^H(kT+v)=\left\{\begin {array}{cc}
\alpha^{(k+2)TH'+n}\hspace{1cm}v+n-1\leqslant T-2\\
\alpha^{(k+3)TH'+n}\hspace{1cm}v+n-1\geqslant T-1\\
\end {array}\right.$$\\
Also by straight calculation from $(3.18)$ we have the same result.$\square$\\

Corresponding to the DT-SIM process, $\{X(\alpha^k), k\in {\bf Z}\}$ with scale
$l=\alpha^T$, $\alpha >1$, $T\in {\bf N }$ there exists a
$T$-dimensional discrete time self-similar Markov process
$W(t)=\big(W^0(t), W^1(t),\ldots,W^{T-1}(t)\big)$ with parameter space
$\check{T}=\{l^{n}; n\in {\bf W}, l=\alpha^T\}$, where \be
W^k(l^n)=W^k(\alpha^{nT})=X(\alpha^{nT+k}),\hspace{1cm}k=0, \ldots,
T-1.\ee
The elements of the covariance matrix which is defined by (3.6) at points $l^n$ and $l^{\tau}$ by (3.10) and (3.11) can be written as\\
$$ Q^H_{jk}(l^n,l^{\tau})=E[W^j(l^{n+\tau})W^k(l^{n})]=\alpha^{2nHT}E[X(\alpha^{\tau T+j})X(\alpha^{k})]$$
\be =\alpha^{2nHT}R_k^H(\tau T+j-k)=\alpha^{2nHT}[\tilde{h}(\alpha^{T-1})]^{\tau}C^H_{jk}R_k^H(0)\ee\\
in which $C^H_{jk}=\tilde{h}(\alpha^{j-1})[\tilde{h}(\alpha^{k-1})]^{-1}$ and $R_k^H(\cdot)$ is defined in (3.10).\\

\begin{theorem}
Let $\{X(\alpha^n),n\in {\bf W}\}$ be a {\em DT-SIM} process with the covariance function $R_n^H(\tau)$.
Also let $\{W(l^n),n\in {\bf W}\}$, defined in $(3.19)$, be its associated $T$-dimensional discrete time self-similar Markov process with covariance function $Q^H(l^n,l^{\tau})$. Then
\be Q^H(l^n,l^{\tau})=\alpha^{2nHT}C_HR_H[\tilde{h}(\alpha^{T-1})]^{\tau},\hspace{7mm}\tau\in{\bf W}\ee
where $\tilde{h}(\cdot)$ is defined by $(3.12)$ and the matrices $C_H$ and $R_H$ are given by $C_H=[C_{jk}^H]_{j,k=0, 1, \ldots, T-1}$, where $C^H_{jk}=\tilde{h}(\alpha^{j-1})[\tilde{h}(\alpha^{k-1})]^{-1}$, and

$$R_H=\left[\begin{array}{cccc}
R^H_0(0)&0&\cdots&0\\
0&R^H_1(0)&\cdots&0\\
\vdots&\vdots&\vdots&\vdots\\
0&0&\cdots&R^H_{T-1}(0)\\
\end{array}\right].\square$$
\end{theorem}

\begin{remark}
It follows from Theorem $3.3$ that for each $k=0, 1, \ldots, T-1$ the process $W^k(l^n)=X(\alpha^{nT+k})$ is a self-similar Markov process for $n\in {\bf W}$. The covariance function of the process is
$$\Gamma^H_k(l^n,l^{\tau})=E[W^k(l^{n+\tau})W^k(l^{n})]=\alpha^{2nHT}[\tilde{h}(\alpha^{T-1})]^{\tau}R^H_k(0),\hspace{1cm}\tau\in {\bf W}$$
where $C^H_{kk}=\tilde{h}(\alpha^{k-1})[\tilde{h}(\alpha^{k-1})]^{-1}=1.$
\end{remark}

The introduced $T$-dimensional self-similar Markov process $W(t)$ with parameter space
$\check{T}=\{l^n, n\in {\bf W}\}$, $l=\alpha^n$, $T\in {\bf N}$, is the counterpart of the
$T$-dimensional stationary Markov process $Y(t)=\big(Y^0(t), Y^1(t), \ldots, Y^{T-1}(t)\big)$.
The spectral density matrix of such $T$-dimensional self-similar process is characterized by the following lemma.
\begin{lemma}
The spectral density matrix ${\bf d}^H(\omega)=[{\bf
d}^H_{jr}(\omega)]_{j,r=0,\ldots,T-1}$ of the $T$-dimensional
self-similar process $\{W(l^n), n\in {\bf W}\}$, defined by $(3.19)$, where $l=\alpha^T$
has the Markov property and is specified by
$${\bf d}^H_{jr}(\omega)=\frac{1}{2\pi}\left[\frac{\tilde{h}(\alpha^{j-1})R_r^H(0)}{\tilde{h}(\alpha^{r-1})
(1-e^{ -i\omega T}\alpha^{- H T}\tilde{h}({\alpha}^{T-1}))}-
\frac{\tilde{h}(\alpha^{r-1})R_j^H(0)}{\tilde{h}(\alpha^{j-1})
\big(1-e^{-i\omega T}\alpha^{HT}\tilde{h}^{-1}({\alpha}^{T-1})\big )}\right]$$
where $R^H_{k}(0)$ is the variance of the process $X(\cdot)$ at point ${\alpha}^{k}$ and $\tilde{h}(\alpha^{k})$ is defined by $(3.12)$.\\
\end{lemma}
{\bf Proof:} As we mentioned prior to the proof of the Theorem 3.1,
we consider $Q^H_{jr}(\cdot,\cdot)$ at discrete points $l^m$ and
$l^s$ where $m,s\in {\bf Z}$, then
$$Q^H_{jr}(l^m,l^s)=E[W^j(l^{m+s})W^r(l^{m})]=l^{2mH}E[W^j(l^s)W^r(1)]=l^{2mH}Q^H_{jr}(l^s).$$
If the $T$-dimensional discrete time self-similar Markov process $W(\cdot)$ is sampled at points $l^n=\alpha^{nT}$, then in (3.7) we have $\tau T$ instead of $\tau$, thus in (3.5) we have $nT$ instead of $n$ and the corresponding spectral density matrix
of the covariance matrix  ${\bf Q}^H( l^s)=[Q_{jr}^H(l^s)]_{j,r=0, 1, \ldots, T-1} $ is
$${\bf d}^H(\omega)=[{\bf d}_{jr}^H(\omega)]_{j,r=0, 1, \ldots, T-1}$$
\vspace{-3mm}
where
$${\bf d}^H_{jr}(\omega)=\frac{1}{2\pi}\sum_{s=0}^{\infty}l^{-Hs}e^{-i\omega sT}Q_{jr}^H(l^s)
+\frac{1}{2\pi}\sum_{s=-\infty}^{-1}l^{-Hs}e^{-i\omega
sT}Q_{jr}^H(l^s):={\bf d}^H_{jr1}(\omega)+{\bf d}^H_{jr2}(\omega).$$
Using (3.20), we evaluate ${\bf d}^H_{jr1}(\omega)$ as
\vspace{-1mm}
$${\bf d}^H_{jr1}(\omega)=\frac{1}{2\pi}\sum_{s=0}^{\infty}e^{-i\omega sT}l^{-  Hs}R_r^H(sT+j-r)
=\frac{1}{2\pi}\sum_{s=0}^{\infty}e^{-i\omega sT}\alpha^{-HsT}[\tilde{h}({\alpha}^{T-1})]^{s}C_{jr}^HR_r^H(0)$$
\be=\frac{\tilde{h}(\alpha^{j-1})R_r^H(0)}{2\pi\tilde{h}(\alpha^{r-1})}\sum_{s=0}^{\infty}
\big(e^{-i\omega T}\alpha^{-HT}\tilde{h}({\alpha}^{T-1})\big)^{s}.\ee
Now we verify the convergence of the above summation. By (3.12) we have
$$|e^{-i\omega T}\tilde{h}(\alpha^{T-1})|=|\tilde{h}(\alpha^{T-1})|=\prod_{j=0}^{T-1}|h(\alpha^j)|
=\prod_{j=0}^{T-1}\left|\frac{R_j^H(1)}{R_j^H(0)}\right|$$
$$=\prod_{j=0}^{T-2}\left|\frac{E[X(\alpha^{j+1})X(\alpha^{j})]}{\sqrt{E[X^2(\alpha^{j+1})]E[X^2(\alpha^j)]}}
\times\frac{E[X(\alpha^{T})X(\alpha^{T-1})]}{\sqrt{E[X^2(1)]E[X^2(\alpha^{T-1})]}}\right|.$$\\
By scale invariance of $X(\cdot)$ we have that $E[X^2(\alpha^{T})]=\alpha^{2TH}E[X^2(1)]$.
Now for $j=0, \ldots,\\ T-1$ if at least one of the $\text{Corr}[X(\alpha^{j+1})X(\alpha^{j})]< 1$ then $\tilde{h}(\alpha^{T-1})<\alpha^{TH}$, and
$$|e^{-i\omega T}\alpha^{- HT}\tilde{h}(\alpha^{T-1})|< 1.$$
Therefore the summation on the right side of (3.22) is convergent. Thus the spectral density is
$${\bf d}^H_{jr1}(\omega)=\frac{\tilde{h}(\alpha^{j-1})R_r^H(0)}{2\pi\tilde{h}(\alpha^{r-1})}
\times\frac{1}{1-e^{-i\omega T}\alpha^{-HT}\tilde{h}({\alpha}^{T-1})}.$$
Now we are to evaluate ${\bf d}^H_{jr2}(\omega)$ as
$${\bf d}^H_{jr2}(\omega)=\frac{1}{2\pi}\sum_{s=1}^{\infty}l^{ Hs}e^{i\omega sT}Q_{jr}^H(l^{-s}).$$
As
$Q_{jr}^H(l^{-s})=E[W^j(l^{-s})W^r(1)]=l^{-2sH}E[W^j(1)W^r(l^{s})]=l^{-2sH}Q_{rj}^H(l^{s})$,
so by a similar method, one can easily verify that
$${\bf d}^H_{jr2}(\omega)=\frac{\tilde{h}(\alpha^{r-1})R_j^H(0)}{2\pi\tilde{h}(\alpha^{j-1})}
\times\frac{e^{ i\omega T}\alpha^{- HT}\tilde{h}({\alpha}^{T-1})}{1-e^{ i\omega T}\alpha^{- HT}\tilde{h}({\alpha}^{T-1})}.$$
So we arrive at the assertion of the lemma.$\square$

\begin{remark}
Lemma $3.4$ provides the spectral density of discrete time self-similar Markov process $\{W^k(l^n), n\in {\bf Z}\}$, defined by $(3.19)$, for $k=0, 1, \ldots, T-1$ as
$${\bf d}^H_{kk}(\omega)=\frac{R_k^H(0)\big(1-\alpha^{-2HT}\tilde{h}^2({\alpha}^{T-1})\big)}{2\pi\big(1-2\cos(\omega T)
\alpha^{- H T}\tilde{h}({\alpha}^{T-1})+\alpha^{-2HT}\tilde{h}^2({\alpha}^{T-1})\big)}.$$
\end{remark}

\begin{remark}
Using Lemma $3.4$, relations $(2.16)$, $(2.17)$ and $(3.20)$, we see that the spectral density matrix $f(\omega)=[f_{jr}(\omega)]_{j,r=0,1,\ldots,T-1}$ of a {\em DT-SIM} process which is defined by $(2.18)$ is fully specified by $\{R_{j}^H(1),R_{j}^H(0),j=0, 1, \ldots, T-1\}$.
\end{remark}

\begin{example}
Here we present the T-dimensional discrete time self-similar Markov process corresponding to the simple Brownian motion, described in Example $3.1$, as $W(l^n)=(W^0(l^n),\\ W^1(l^n), \ldots, W^{T-1}(l^n))$, where $W^k(l^n)=X(\alpha^{nT+k})$. Now as we mentioned in Lemma $3.5$ we obtain spectral density matrix of $W(l^n)$. In Example $3.1$ we find that $\tilde{h}(\alpha^{j-1})=1$ and $\tilde{h}(\alpha^{r-1})=1$ as $j, r= 0, 1, \ldots, T-1$, $\tilde{h}(\alpha^{T-1})=\alpha^{TH'}$, $H'=H-\frac{1}{2}$ and $R^H_r(0)=\alpha^{2TH'+r}$, $R^H_j(0)=\alpha^{2TH'+j}$ thus the spectral density matrix of $W(l^n)$ is
$${\bf d}^H_{jr}(\omega)=\frac{\alpha^{2TH'}}{2\pi}\left[
\frac{\alpha^{r}}{1-e^{-i\omega T}\alpha^{-T/2}}-\frac{\alpha^{j}}{1-e^{-i\omega T}\alpha^{T/2}}
\right].$$
\end{example}

\hspace{-6mm}{\Large\bf Acknowledgements}

The authors would like to express their thanks to both anonymous referees for valuable comments and suggestions which improve the original manuscript.

\end{document}